\documentclass[12pt]{article}

      \usepackage{latexsym}
         \usepackage[reqno, namelimits, sumlimits]{amsmath}
         \usepackage{amssymb, amsfonts}
         \usepackage{amsthm}

 \newtheorem{theorem}{Theorem}[section]

 \newtheorem{pro}[theorem]{Proposition}

\title{Time decay for solutions to the Stokes equations with drift
}

 \author{M.~Schonbek, G.~Seregin}

\date{}
\begin{document}
\maketitle
\begin{abstract}
In this note, we study the behaviour of Lebesgue norms $\|v(\cdot,t)\|_p$ of solutions $v$ to the Cauchy problem for the Stokes system with drift $u$,  which is supposed to be a divergence free smooth vector valued function  satisfying a scale invariant condition.
 
 \end{abstract}

\setcounter{equation}{0}
\section{Inroduction}

The main aim of the paper is the following Stokes system with a drift $u$

\begin{equation}\label{moreperturbproblem}
    \partial_tv-u\cdot \nabla v-\triangle v-\nabla q=-{\rm div}\, F,\qquad {\rm div}\,v=0
\end{equation}
in $Q_+=\mathbb R^3\times ]0,\infty[$ and
\begin{equation}\label{moreinitialldata}
    v(x,0)=0
\end{equation}
for $x\in \mathbb R^3$.

It is supposed that a tensor-valued field  $F$ is smooth and compactly supported in $Q_+$. In addition, let us assume that $F$ is skew symmetric and therefore
\begin{equation}
	\label{compatibility2}
	{\rm div}\,{\rm div}\,F=0.
\end{equation} 

As to the drift $u$, one may assume that $u$ is a bounded divergence free field in $Q_+$, say $|u|\leq 1$ there,  whose derivatives of any order exist and are bounded in $Q_+$. 
 
 It is not so difficult to prove, see Appendix I, the following statement. 
\begin{pro}\label{2p1} There exists a unique solution $v$ to (\ref{moreperturbproblem}) and (\ref{moreinitialldata}) with properties:
$$\nabla^l\partial^k_t v\in L_2(Q_+)$$
for $k,l=0,1,...$ except $k+l=0$,
$$\nabla^{l+1}\partial^k_t q\in L_2(Q_+)$$
for $k,l=0,1,...$,
$$v\in L_{2,\infty}(Q_+), \qquad q\in L_{2,\infty}(Q_+)$$
for any $k=0,1,...$.
\end{pro} 

The goal of the paper is to study how $L_p$-norms of the velocity field $v$ $\Big(\|v(\cdot,t)\|_p:=\Big(\int\limits_{\mathbb R^3}|v(x,t)|^pdx\Big)^\frac 1p)$ behave as $t\to\infty$. In particular, two cases are of great interest: $p=1$ and $p=2$. 
 
 Let us impose a decay assumption on the drift
 \begin{equation}
 	\label{driftdecay}
 	|u(x,t)|\leq \frac {c_d}{|x|+\sqrt {t}}
 \end{equation}
 for all $(x,t)\in Q_+$.
 
 Two results will be proven in the paper.
 \begin{theorem}
 	\label{Time Decay} Let $v$ be 
 	a  solution $v$ to (\ref{moreperturbproblem}) and (\ref{moreinitialldata}) and let $u$ satisfy (\ref{driftdecay}). Then for 
 any $m=0,1...$,  two decay estimates are valid: 
\begin{equation}\label{L1}
  \|v(\cdot,t)\|_1\leq c(m,c_d)\sqrt{t}^\frac 32\frac 1{\ln^m(t+e)}
\end{equation}
and
\begin{equation}\label{L2}
  \|v(\cdot,t)\|_2\leq \frac {c(m,c_d)}{\ln^m(t+e)}.
\end{equation}	
 	
 \end{theorem}

 To motivate the aforesaid problem and the assumptions made, consider the Navier-Stokes system 
 $$\partial_tw+w\cdot\nabla w-\Delta w=-\nabla r,\qquad {\rm div}\,w=0$$
 in the unit  parabolic ball $Q=B\times ]-1,0[$ for 
 functions 
 $w\in L_\infty(-1,0;L_2(B))\cap L_2(-1,0;W^1_2(B))$ and $r
 \in L_{\frac 32}(Q)$ satisfying the additional restriction 
 \begin{equation}
 	\label{typeI}
 	|w(x,t)|\leq \frac {c_d}{|x|+\sqrt{-t}}
 \end{equation}
 for all $(x,t)\in Q$. Our aim is to understand whether or not the origin $z=(x,t)=(0,0)$ is a regular point of $w$,
 i.e., there exists $\delta>0$ such that $v$ is essentially bounded  in the parabolic ball $Q(\delta)=B(\delta)\times ]-\delta^2,0[$. Here, as usual, $B(r)$ stands for the ball of radius $r$ centered at the origin.
 The answer  is certainly positive if $c_d$ is sufficiently small. However, we would not like to make such an assumption at this point. In \cite{SerSve09},  it has been shown that if $z=0$ is a singular point of $w$ then a so-called a mild bounded ancient solution $\tilde u$ to the Navier-Stokes equations in $Q_-=\mathbb R^3\times ]-\infty,0[$ exists and it is non-trivial. The latter means the following: $\tilde u\in L_\infty(Q_-)$ ($|\tilde u|\leq 1$ a.e. in $Q_-$and $|u(0)|=1$) and there exists a scalar function $\tilde p\in L_\infty(-\infty,0;BMO(\mathbb R^3))$ such that the pair $\tilde u$ and $\tilde p$ satisfy the classical Navier-Stokes system
\begin{equation}\label{NSS}
    \partial_t \tilde u+\tilde u\cdot\nabla\tilde  u-\triangle \tilde u=-\nabla \tilde p,\qquad {\rm div}\,\tilde u=0
\end{equation}
in $Q_-$ in the sense of distributions. It is known, see \cite{KNSS2009}, that $\tilde u$ is infinitely smooth and all its derivatives are bounded.  Moreover, it can be shown, see Appendix II, that, for $u(x,t)=\tilde u(x,-t)$,
\begin{equation}
	\label{main formula}
	\int\limits_{Q_+} u\cdot  {\rm div}\,F dxdt=-\lim\limits_{T\to \infty}\int\limits_{\mathbb R^3} u(x,T)\cdot v(x,T)dx.
\end{equation}
If time decay of $v$ is such that, for any tensor-valued field  $F
\in C^\infty_0(\mathbb R^3)$,
 obeying condition (\ref{compatibility2}),
  the limit on the right hand side of (\ref{main formula}) vanishes, then one can easily show that $u$ must be a function of time only. Indeed, 
we then have $$\int\limits_{Q_+}\nabla u:Fdxdt=0.$$
The latter means that the skew symmetric part of $\nabla u$  vanishes in $Q_+$. Since $u$ is a divergence free field,  $u$ is a bounded  harmonic function and so does $\tilde u$ in $Q_-$. But $\tilde u$ is a bounded mild ancient solution to the Navier-Stokes equation and thus  must be a constant in $Q_-$ as well as $u$ in $Q_+$. But condition (\ref{driftdecay}) means that $\tilde u$ is identically zero. This finally would prove that $z=0$ is not a singular point of $w$ and condition (\ref{typeI}) is in fact a regularity condition.

 Unfortunately, decay bounds in Theorem \ref{Time Decay} do not provide the above scenario. 
Let us give a couple of  bounds on $c_d$ that give
 a required time decay.

 To describe  the first case, we are going to use a solution formula for the Stokes system with non-divergence free right hand side. 
 
 Let
$$\mathcal F=-v\otimes u+F.$$
The solution to problem (\ref{moreperturbproblem}), (\ref{moreinitialldata}) has the form, see for instance  \cite{KNSS2009}, 
\begin{equation}\label{solutionformulaII}
  v(x,t)=\int\limits^t_0\int\limits_{\mathbb R^3}K(x-y,t-s)\mathcal F(y,s)dy ds,
\end{equation}
where the potential $K=(K_{ijl})$ defined with the help of the standard heat kernel in the following way
$$\Delta \Phi(x,t)=\Gamma(x,t)$$
and
$$K_{ijl}=\Phi_{,ijl}-\delta_{il}\Phi_{,kkj}.$$
It is easy to check that the following bound is valid:
\begin{equation}\label{kernelboundII}
 | K(x,t)|\leq \frac {c_1}{(t+|x|^2)^2}
\end{equation}
and therefore
\begin{equation}\label{intkerboundII}
  \int\limits_{\mathbb R^3}|K(x,t)|dx\leq \frac {c_*}{\sqrt {t}}
\end{equation}
 with $c_*=cc_1$, where $c$ is an absolute constant. 
 \begin{theorem} 
 \label{smallnessM}
 Assume that 
 \begin{equation}
 	\label{smallM}
 	4c_*c_d<1.
 \end{equation}
 Then 	
\begin{equation}
	\label{innerto0}
\int\limits_{\mathbb R^3}v(x,T)\cdot u(x,T)dx\to0\end{equation}
as $T\to\infty$.
 \end{theorem}

 To descibe the second case,  let us  introduce the operator
 $K:\mathcal L_2
 \to J_2
 $, where
 $\mathcal L_2
 $ consists of all tensor-valued functions, belonging to  $L_2(\mathbb R^3)$ and satisfying condition (\ref{compatibility2}), and $J_2$ is a space of square integrable  divergence free fields in $\mathbb R^3$. The action of this operator is defined as $A_F=KF$, where $A_F$ is the unique solution to the 
 following problem
 \begin{equation}
 	\label{elliptic stuff1}
 	{\rm rot}\, A_F=-{\rm div}\,F. 
 \end{equation}
 The elliptic theory reads that operator $K$  is bounded. 

In addition, one may introduce the second operator
$M:L_2(\mathbb R^3;\mathbb M^{3\times3})\to L_2(\mathbb R^3)$  so that 
\begin{equation}
	\label{elliptic stuff2}
\Delta q_F=-{\rm div}\,{\rm div}\,F,\end{equation}
where
$q_F=MF$. 

Actually, we have fixed the pressure  $q=q_{v\otimes u}$ in Proposition \ref{2p1}. This will be done everywhere in what follows. Our result is the following.
\begin{theorem}
\label{smallness}
Let 
$$c_d\leq \frac {\sqrt 3}{2\|K\|(1+\sqrt 3\|M\|)}.$$	Then (\ref{innerto0}) is true.
\end{theorem}


\setcounter{equation}{0}
\section{Time Decay of $L_1$-Norm}

Now,  from (\ref{solutionformulaII}), it follows
$$\|v(\cdot,t)\|_p\leq\int\limits^t_0\|\int\limits_{\mathbb R^3}|K(\cdot-y,t-s)|\mathcal F(y,s)dy\|_p ds
$$
Applying H\"older inequality and taking into account  (\ref{intkerboundII}), we find
$$\|v(\cdot,t)\|_p\leq\int\limits^t_0\Big(\int\limits_{\mathbb R^3}\Big(\int\limits_{\mathbb R^3}K(x-y',t-s)dy'\Big)^\frac p{p'}\times$$$$\times\int\limits_{\mathbb R^3}|K(x-y,t-s)||\mathcal F(y,s)|^pdydx\Big)^{\frac 1{p}} ds\leq c\int\limits^t_0\frac 1{\sqrt{t-s}}\|\mathcal F(\cdot,s)\|_pds
$$
for any $p\geq 1$.

Now, for $p$,
 satisfying the condition
\begin{equation}\label{rangefor pII}
  p\in ]6/5,2[,
\end{equation}
 H\"older inequality gives the following estimate
$$\|v(\cdot,t)\|_{p}\leq c\int\limits^t_0\frac {ds}{\sqrt{t-s}}
\Big(\int\limits_{\mathbb R^3}|\mathcal F(y,s)|^2(\sqrt{s}+|y|)^2dy\Big)^\frac 12\times $$$$\times\Big(\int\limits_{\mathbb R^3}\Big(\frac 1{\sqrt{s}+|y|}\Big)^{\frac {2p}{2-p}}dy\Big)^\frac {2-p}{2p}.
$$
By changing variables $y=z\sqrt{s}$, 
$$\Big(\int\limits_{\mathbb R^3}\Big(\frac 1{\sqrt{s}+|y|}\Big)^{\frac {2p}{2-p}}dy\Big)^\frac {2-p}{2p}\leq
\sqrt{s}^{-\frac{5p-6}{2}}\Big(\int\limits_{\mathbb R^3}\Big(\frac 1{1+|z|}\Big)^{\frac {2p}{2-p}}dz\Big)^\frac {2-p}{2p}= $$
$$=  C(p)\sqrt{s}^{-\frac{5p-6}{2p}}$$
with
$$C(p):=\Big(\int\limits_{\mathbb R^3}\Big(\frac 1{1+|z|}\Big)^{\frac {2p}{2-p}}dz\Big)^\frac {2-p}{2p}\to \infty $$
as $p\to 6/5+0$. So,
$$\|v(\cdot,t)\|_{p}\leq C(p)\int\limits^t_0\frac {ds}{\sqrt{t-s}}\sqrt{s}^{-\frac{5p-6}{2p}}
\Big(\int\limits_{\mathbb R^3}|\mathcal F(y,s)|^2(\sqrt{s}+|y|)^2dy\Big)^\frac 12.
$$

Now, by our assumptions on $F$ and by (\ref{driftdecay}), 
\begin{equation}\label{bound}
  \int\limits_{\mathbb R^3}|\mathcal F(y,s)|^2(\sqrt{s}+|y|)^2dy\leq c(c_d\|v(\cdot,s)\|_2+\|G(\cdot,s)\|_{2})^2,
  \end{equation}
where $
G(y,s)=F(y,s)(\sqrt{s}+|y|),
$
and thus
$$\|v(\cdot,t)\|_{p}\leq C(p)\int\limits^t_0\frac {ds}{\sqrt{t-s}}\sqrt{s}^{-\frac{5p-6}{2p}}(c_d\|v(\cdot,s)\|_2+\|G(\cdot,s)\|_{2})$$
$$\leq C(p)A_p(s)$$
with
\begin{equation}\label{formulaforA}
  A_p(s):=\int\limits^t_0\frac {ds}{\sqrt{t-s}}\sqrt{s}^{-\frac{5p-6}{2p}}(c_d\|v(\cdot,s)\|_2+\|G(\cdot,s)\|_{2}).
\end{equation}
So,
\begin{equation}\label{intermediateII}
  \|v(\cdot,t)\|_{p}\leq C(p)A_p(t).
\end{equation}

Now, one can repeat the above arguments  for  $p=1$ and find
$$\|v(\cdot,t)\|_{1}\leq \int\limits^t_0\frac c{\sqrt{t-s}}\int\limits_{\mathbb R^3}|\mathcal F(y,s)|dy ds.$$
Since
$$|\mathcal F(y,s)|\leq c\frac {c_d|v(y,s)|+|G(y,s)|}{\sqrt{s}+|y|},$$
the latter estimate can be transform as follows:
$$\|v(\cdot,t)\|_{1}\leq c\int\limits^t_0\frac {ds}{\sqrt{t-s}}\int\limits_{\mathbb R^3}\frac {c_d|v(y,s)|+|G(y,s)|}{\sqrt{s}+|y|}dy\leq $$
$$\leq c\int\limits^t_0\frac {ds}{\sqrt{t-s}}\Big(\int\limits_{\mathbb R^3}\Big(\frac 1{\sqrt{s}+|y|}\Big)^\frac {6+5\varepsilon}{1+5\varepsilon}dy\Big)^\frac {1+5\varepsilon}{6+5\varepsilon} \Big(\int\limits_{\mathbb R^3}(c_d|v(y,s)|+$$$$+|G(y,s)|)^\frac {6+5\varepsilon}5 dy\Big)^\frac 5{6+5\varepsilon}$$
for some positive $0<\varepsilon<3/10$.
Hence,
$$\|v(\cdot,t)\|_{1}\leq C_1(\varepsilon)\int\limits^t_0\frac {ds}{\sqrt{t-s}}\sqrt{s}^{ 3\frac{1+5\varepsilon}{6+5\varepsilon}-1}\Big(\int\limits_{\mathbb R^3}(c_d|v(y,s)|+$$$$+|G(y,s)|)^\frac {6+5\varepsilon}5 dy\Big)^\frac 5{6+5\varepsilon}$$
with $$C_1(\varepsilon):=\Big(\int\limits_{\mathbb R^3}\Big(\frac 1{1+|z|}\Big)^\frac {6+5\varepsilon}{1+5\varepsilon}dz\Big)^\frac {1+5\varepsilon}{6+5\varepsilon}.$$
Simplifying slightly the previous bound, we have
$$\|v(\cdot,t)\|_{1}\leq C_1(\varepsilon)\int\limits^t_0\frac {ds}{\sqrt{t-s}}\sqrt{s}^{\frac {-3+10\varepsilon}{6+5\varepsilon}}(\|c_dv(\cdot,s)\|_{\frac {6+5 \varepsilon}5}+\|G(\cdot,s)\|_\frac {6+5\varepsilon}5)dy.$$
By (\ref{intermediateII}),
$$\|v(\cdot,s)\|_{\frac {6+5 \varepsilon}5}\leq C( 6/5+\varepsilon) A_{\frac 65+\varepsilon}(t).$$
So, the final estimate of $L_1$-norm is:
$$\|v(\cdot,t)\|_{1}\leq C_3(\varepsilon)\int\limits^t_0\frac {ds}{\sqrt{t-s}}\sqrt{s}^{\frac {-3+10\varepsilon}{6+5\varepsilon}}(c_d A_{\frac 65+\varepsilon}(s)+$$\begin{equation}\label{L1-norm}+\|G(\cdot,s)\|_\frac {6+5\varepsilon}5)dy\end{equation}
 with $C_3(\varepsilon)\to\infty$ as $\varepsilon\to 0$.
 
Since the energy of $v$ is bounded,  one  can derive
 from (\ref{formulaforA}) the following:
$$A_p(t)\leq c(s)(c_d\|v\|_{2,\infty}+\|G\|_{2,\infty})\sqrt{t}^\frac {6-3p}{2p}$$
and thus
$$A_{\frac 65+\varepsilon}(t)\leq c(\varepsilon)(c_d\|v\|_{2,\infty}+\|G\|_{2,\infty})\sqrt{t}^\frac{12-15\varepsilon}{2(6+5\varepsilon)}.$$
Now, (\ref{L1-norm}) is giving to us:
$$\|v(\cdot,t)\|_1\leq c\sqrt{t}^\frac 32$$
where $c$ depends on the data of the problem.

\setcounter{equation}{0}
\section{Improvement for  $L_2$-norm}

Here, we are going to use methods developed in 
\cite{Schonbek1985} and \cite{Wiegner1987}, see also \cite{BorMia1988} and \cite{KajMiy1986}.

We have the energy inequality
\begin{equation}\label{energyineq}
  \partial_ty(t)+\|\nabla v(\cdot,t)\|^2_2\leq \|F(\cdot,t)\|^2_2
\end{equation}
with $y(t)=\|v(\cdot,t)\|^2_2$.

The Fourier transform and Plancherel identity give us
$$ \partial_ty(t)\leq -\int\limits_{\mathbb R^3}|\xi|^2|\widehat{v}(\xi,t)|^2d\xi+ \|F(\cdot,t\|^2_2\leq $$
$$=-\int\limits_{|\xi|>g(t)}|\xi|^2|\widehat{v}(\xi,t)|^2d\xi-\int\limits_{|\xi|\leq g(t)}|\xi|^2|\widehat{v}(\xi,t)|^2d\xi+ \|F(\cdot,t\|^2_2,$$
where $g(t)$ is a given function which will be specified later on.
The latter implies
$$y'(t)+g^2(t)y(t)\leq \int\limits_{|\xi|\leq g(t)}(g^2(t)-|\xi|^2)|\widehat{v}(\xi,t)|^2d\xi+ \|F(\cdot,t\|^2_2.$$

Taking the Fourier transform of the Navier-Stokes equation, we find
$$\partial_t\widehat{v}+|\xi|^2\widehat{v}=-\widehat{H},$$
where
$$H=-\mbox{div}\,(v\otimes u+\mathbb Iq-F).$$
Clearly,
$$\widehat{v}(\xi,t)=-\int^t_0\exp\{-|\xi|^2(t-s)\}\widehat{H}(\xi,s) ds$$
and
$$|\widehat{H}(\xi,s)|\leq |\xi|\Big(\||v(\cdot,s)||u(\cdot,s)|\|_1+\|F(\cdot,s)\|_1\Big).$$

Denoting
$$k(t)=\||v(\cdot,t)|\|_1,$$
we notice
$$\||v(\cdot,s)||u(\cdot,s)|\|_1\leq \sqrt{s}^{-1}c_dk(s).$$
So,
$$|\widehat{v}(\xi,t)|\leq c\int^t_0\exp\{-|\xi|^2(t-s)\}|\xi|(\sqrt{s}^{-1}k(s)+\|F(\cdot,s)\|_1)ds.$$
Applying the H\"older inequality, we show
$$y'(t)+g^2(t)y(t)\leq $$$$\leq c\int\limits_{|\xi|\leq g(t)}(g^2(t)-|\xi|^2)\Big(\int^t_0\exp\{-|\xi|^2(t-s)\}|\xi|(\sqrt{s}^{-1}k(s)+\|F(\cdot,s\|_1)ds\Big)^2\leq $$
$$\leq c\int^t_0(s^{-1}k^2(s)+\|F(\cdot,s\|^2_1)ds\times$$
$$\times\int^t_0\int\limits_{|\xi|\leq g(t)}(g^2(t)-|\xi|^2)|\xi|^2\exp\{-2|\xi|^2(t-s_1)\}d\xi ds_1+\|F(\cdot,t)\|_2^2.$$
The latter integral can be estimated in the following way:
$$\int\limits^t_0\int\limits_{|\xi|\leq g(t)}(g^2(t)-|\xi|^2)|\xi|^2\exp\{-|\xi|^2(t-s_1)\}d\xi ds_1=$$
$$=c\int\limits^t_0\int\limits^{g(t)}_0(g^2(t)-r^2)r^4\exp\{-2r^2(t-s_1)\}dr ds_1\leq
$$$$\leq cg^6(t)\int\limits^t_0\int\limits^{g(t)}_0\exp\{-2r^2(t-s_1)\} d(r\sqrt{t-s_1})\frac {ds_1}{\sqrt{t-s_1}}\leq $$
$$\leq cg^6(t)\int\limits^t_0
\frac {ds_1}{\sqrt{t-s_1}}\int\limits^\infty_0\exp\{-2z^2\}dz\leq cg^6(t)\sqrt{t}.$$

Coming back to our energy inequality, we find
$$y'(t)+g^2(t)y(t)\leq $$$$\leq K(t):=cg^6(t)\sqrt{t}\int^t_0(s^{-1}k^2(s)+\|F(\cdot,s\|^2_1)ds+\|F(\cdot,t)\|_2^2.$$
Then Grownwall inequality implies
$$
y(t)\leq c\int\limits^t_0\exp\Big\{-\int\limits^t_sg^2(\tau)d\tau \Big\}K(s)ds.$$

\setcounter{equation}{0}
\section{Proof  of Theorem \ref{Time Decay}}

The proof is on induction in $m$. The basis of induction has been already established in Section II.
Let us assume that our statement is true for $m$ and show that it is true for $m+1$.

We can  estimate $K(t)$ using the  fact that $F$ has a compact support
$$K(t)\leq cg^6(t)\sqrt{t}\int^t_0(\sqrt{s}\ln^{-2m}(s+e)+\|F(\cdot,s)\|^2_1)ds+\|F(\cdot,t)\|_2^2\leq $$

$$\leq C(\|F\|_{1,\infty},m)g^6(t)\sqrt{t}\int^t_0\sqrt{s}\ln^{-2m}(s+e)ds+\|F(\cdot,t)\|_2^2\leq$$
$$\leq C(\|F\|_{1,\infty},m)g^6(t)t^2\ln^{-2m}(t+e)+\|F(\cdot,t)\|_2^2.$$

Let
\begin{equation}\label{functiong}
  g^2(t)=\frac {h'(t)}{h(t)}.
\end{equation}
Then
$$\int\limits^t_0\exp\Big\{-\int\limits^t_sg^2(\tau)d\tau \Big\}(g^6(s)s^2\ln^{-2m}(s+e)+\|F(\cdot,s)\|_2^2)ds=$$
$$=\frac 1{h(t)}\int\limits^t_0\Big(\frac {s^2\ln^{-2m}(s+e)}{h^2(s)}(h'(s))^3+h(s)\|F(\cdot,s)\|_2^2\Big)ds.$$
Now,
one specify function $g$ by a particular choice of fuction h, setting
\begin{equation}\label{functionhwithk}
  h(t)=\ln^k(t+e)
\end{equation}
 for some $k>2m+2$. Then
$$\frac 1{h(t)}\int\limits^t_0\frac {s^2\ln^{-2m}(s+e)}{h^2(s)}(h'(s))^3ds=$$$$ =\frac 1{\ln^k(t+e)}
\int\limits^t_0\frac {s^2\ln^{-2m}(s+e)}{(s+e)^3}k^3\ln^{k-3}(s+e)ds\leq$$
$$ =\frac 1{\ln^k(t+e)}
\int\limits^t_0\frac {s^2}{(s+e)^3}k^3\ln^{k-2m-3}(s+e)ds\leq$$
$$\leq \frac 1{\ln^k(t+e)}\frac {k^3}{k-2m-2}(\ln^{k-2m-2}(t+e)-1)\leq $$$$\leq c(k,m)\frac 1{\ln^{2m+2}(t+e)}.
$$
Since $s\mapsto \|F(\cdot,s)\|_2^2$ has a compact support in $]0,\infty[$,
we find
$$\|v(\cdot,t)\|_2\leq \frac c{\ln^{m+1}(t+e)}.$$
Then, as it follows from (\ref{formulaforA}),
$$A_p(t)\leq C(\|G\|_{2,\infty},p,m)\sqrt{t}^\frac {6-3p}{2p}\frac 1{\ln^{m+1}(t+e)}$$
and thus
$$A_{\frac 65+\varepsilon}(t)\leq C(\|G\|_{2,\infty},\varepsilon,m)\sqrt{t}^\frac{12-15\varepsilon}{2(6+5\varepsilon)}\frac 1{\ln^{m+1}(t+e)}.$$
And again from (\ref{L1-norm}), it follows finally that
$$\|v(\cdot,t)\|_1\leq c\sqrt{t}^\frac 32\frac 1{\ln^{m+1}(t+e)}.$$

\setcounter{equation}{0}
\section{Liouville type theorems} 
\textsc{Proof of Theorem \ref{smallnessM}} From (\ref{solutionformulaII}) and from (\ref{driftdecay}), one can derive
\begin{equation}
\label{mainiteration}	
f(t)\leq c_*\int\limits^t_0\frac 1{\sqrt{t-s}}\Big(\frac {c_d}{\sqrt{s}}f(s)+\|F(\cdot,s)\|_1\Big)ds,
\end{equation}
where $f(t):=\|v(\cdot,t)\|_1$. Since $F$ is compactly supported, (\ref{mainiteration}) can be reduced to  the following form:
$$f(t)\leq A+c_*c_d\int\limits_0^t\frac 1{\sqrt{t-s}\sqrt{s}}f(s)ds.$$
Now, fix an arbitrary $T>0$. Then, for any $t\in ]0,T]$, we have
$$f(t)\leq A+4c_*c_dM(T),$$
where
$M(T)=\sup\limits_{0<t\leq T}f(t)$. Hence,
$$M(T)\leq A+4c_*c_dM(T)$$
for any $T>0$. Finally, we see that
$$\|v(\cdot,t)\|_1\leq c=\frac A{1-4c_*c_d}$$
for all $t>0$.
Therefore,
$$\Big|\int\limits_{\mathbb R^3}u(x,t)\cdot v(x,t)dx|\leq \frac c{\sqrt{t}}\to0$$
as $t\to\infty$. $\Box$


\textsc{Proof of Theorem \ref{smallness}}  Assume that $F$ is skew symmetric and therefore satisfies condition (\ref{compatibility2}).

 
 Equation (\ref{moreperturbproblem}) can be written as follows: 
 \begin{equation}
 \label{newform3}
 \partial_tv-\Delta v={\rm div}\,F_0,	
 \end{equation}
where 
$$F_0=v\otimes u+\nabla q\mathbb I-F.$$
We know from previous results
that
\begin{equation}
	\label{Fo}
	F_0\in L_{2,\infty}(Q_+),\qquad {\rm div}\,F_0\in L_2(Q_+).
\end{equation}

Since ${\rm div}\,{\rm div\,}F_0=0$, we can apply the elliptic theory and conclude that there exists a divergence free field $A(\cdot,t)$
such that 
\begin{equation}
\label{divrot}
{\rm rot}\,A(\cdot,t)={\rm div}\,F_0(\cdot,t)	
\end{equation}
in $\mathbb R^3$
and the following estimate holds
\begin{equation}
	\label{estdivrot}
	\|A(\cdot,t)\|_2\leq \|K\|\|F_0(\cdot,t)\|_2
\end{equation}
for all $t\in ]0,\infty[$.
Taking into account the definition of the operator $M$,
one can go further and derive from (\ref{estdivrot})
$$	\|A(\cdot,t)\|_2\leq \|K\|(\|v(\cdot,t)\otimes u(\cdot,t)\|_2+\|q_{v\otimes u}(\cdot,t)\mathbb I\|_2+\|F(\cdot,t)\|_2)\leq
$$
$$\leq \|K\|(1+\sqrt{3}\|M\|)
\|v(\cdot,t)\otimes u(\cdot,t)\|_2+h(t),$$
where $h(t)=\|K\|\|F(\cdot,t)\|_2$
and thus
\begin{equation}
	\label{estdivrot1}
	\|A(\cdot,t)\|_2\leq \|K\|(1+\sqrt{3}\|M\|)
	\frac {c_d}{\sqrt{t}}
\|v(\cdot,t)\|_2+h(t)
\end{equation}


With the above $A$, let us consider the 
Cauchy problem 
\begin{equation}
	\label{bproblem}
	\partial_tB-\Delta B=A
\end{equation}
\begin{equation}
	\label{binit}
	B(\cdot,0)=0.
\end{equation}
Problem (\ref{bproblem}), (\ref{binit}) has a unique solution defined for all positive $t$ and $B\in W^{2,1}_2(Q_T)$ for all $T>0.$ Since $A(\cdot,t)$ is divergence free, so is $B(\cdot,t)$. Now, let $w={\rm rot}\, B$. Then we can see that $w$ is a  solution to equation (\ref{newform3}) and sinse it vanishes at $t=0$, we can state that $w=v$.

Now, let us analyse the Cauchy problem for $B$. It is easy to see that $B$ satisfies the energy identity
\begin{equation}
	\label{benergy}
	\frac 12 \partial_t\|B(\cdot,t)\|^2_2+\|\nabla B(\cdot,t)\|^2_2=\int\limits_{\mathbb R^3}A(x,t)\cdot B(x,t)dx.
	\end{equation}
	Taking into account the simple identity
	$$
	\|v(\cdot,t)\|_2=\|\nabla B(\cdot,t)\|_2, 
	$$ one can derive from (\ref{estdivrot1}) the following estimate
	$$  	\frac 12 \partial_t\|B(\cdot,t)\|^2_2+\|v(\cdot,t)\|_2^2   \leq \|K\|(1+\sqrt{3}\|M\|)
	\frac {c_d}{\sqrt{t}}
\|v(\cdot,t)\|_2\|B(\cdot,t)\|_2+$$$$+h(t)\|B(\cdot,t)\|_2.$$
Applying the Young inequality, we find 

$$	\frac 12 \partial_t\|B(\cdot,t)\|^2_2\leq
\|K\|^2(1+\sqrt{3}\|M\|)^2\frac {c_d^2}{4t}\|B(\cdot,t)\|_2^2+\frac 12 h(t)(\|B(\cdot,t)\|_2^2+1)
$$
Let us introduce the important constant
$$l=     \|K\|^2(1+\sqrt{3}\|M\|)^2\frac {c_d^2}{2}.$$
Then the previous inequality leads to
$$\|B(\cdot,t)\|^2_2\leq t^l\int\limits^t_0\frac {h(\tau)}{\tau^l}\exp{\Big(-\int^t_\tau h(s)ds}\Big)d\tau.
$$
Taking into account that $F$ is compactly supported in $Q_+$, we have
$$\|B(\cdot,t)\|_2^2\leq c_Ft^l.$$
From here, it is easy to derive the following:
\begin{equation}
	\label{Sohr}
	\int\limits_o^t\|v(\cdot,s)\|^2_2ds\leq c_Ft^l.
\end{equation}
We denote all the constant depending of $F$ and its support by $c_F$.

Having estimate (\ref{Sohr}) in mind, let us go back to equation (\ref{newform3}) multiplying it by $tv$ and integrating result over $\mathbb R^3\times ]0,t[$, as a result, we find the following differential inequality
$$\frac 12t\|v(\cdot,t)\|^2_2+\int\limits^t_0\|\nabla v(\cdot,s)\|^2_2ds =\frac 12 \|v(\cdot,t)\|^2_2+\int\limits^t_0\int\limits_{\mathbb R^3}sF(x,s)\cdot v(x,s) ds\leq$$
$$\leq c_F(\int\limits^t_0\|v(\cdot,s)\|_2^2ds+1).$$
The latter, together with boundedness of $\|v(\cdot,t)\|_2$,
implies the bound
$$\|v(\cdot,t)\|^2_2\leq c_F(t+1)^{l-1},$$
which, in turn, allows to improve the decay of $\|v(\cdot,t)\|_1$.
To this end, we are going back to (\ref{intermediateII}) and (\ref{L1-norm}). Indeed, by the assumption of the theorem $l<3/4$,
$$A_p(t)\leq c\int\limits^t_0\frac 1{\sqrt{t-s}}s^{-\frac {5p-6}{4p}}(s+1)^{l-1}ds\leq c\int\limits^t_0\frac 1{\sqrt{t-s}}s^{-\frac {5p-6}{4p}+l-1} ds\leq  $$
$$\leq c t^{\frac {6-3p}{4p} +l-1}.$$
Letting $p=6/5+\varepsilon$, for sufficiently small positive $\varepsilon$, we find
$$\|v(\cdot,t)\|_1\leq c(\sqrt{t})^{\frac 32+2(l-1)}.$$
This shows 
$$\Big|\int\limits_{\mathbb R^3}v(\cdot,t)\cdot u(\cdot,t)dx\Big|\leq c(\sqrt{t})^{\frac 12+2(l-1)}\to0$$ as $t\to\infty$ provided $l<\frac 34$. $\Box$

\setcounter{equation}{0}
\section{Appendix I}

\textsc{Proof} We recall that all derivatives of $u$ are bounded.

First of all, there exists a unique energy solution. This follows from the identity
$$\int\limits_{Q_+}(u\cdot\nabla v)\cdot vdx dt=0$$
and from the inequality
$$\Big|-\int\limits_{Q_+}{\rm div}\,F\cdot vdx dt\Big|=\Big|\int\limits_{Q_+}F:\nabla v dx dt\Big|\leq \Big(\int\limits_{Q_+}|F|^2dx dt\Big)^\frac 12\Big(\int\limits_{Q_+}|\nabla v|^2dx dt\Big)^\frac 12$$
So, we can state that
\begin{equation}\label{energyestimate}
    v\in L_{2,\infty}(Q_+),\qquad \nabla v\in L_2(Q_+).
\end{equation}
The latter means that $u\cdot \nabla v\in L_2(Q_+)$.
The pressure can be recovered from the pressure equation 
$$\triangle q={\rm div}\,{\rm div}\,(F-v\otimes u).$$ 
One of solutions to the above equation has the form
$$q_0(x,t)=-\frac 13v(x,t)\cdot u(x,t)+\lim\limits_{\varepsilon\to0} \int\limits_{|x-y|>\varepsilon}\nabla^2E(x-y):v(y,t)\otimes u(y,t)dy,$$ where $E$
is the fundamental solution to the Laplace operator. All others differs from $q_0$ by a function of time only. Let us fix the pressure by setting $q=q_0$. The theory of singular integrals implies that
$$q\in L_{2,\infty}(Q_+), \qquad \nabla q\in L_2(Q_+). $$
Then, by  properties of solutions to the heat equation, we have
$$\nabla^2 v\in  L_2(Q_+)\qquad  \partial_t v\in L_2(Q_+).$$
Going back to the pressure equation, let us re-write it in the following way
$$\triangle q={\rm div}\,{\rm div}\,F -u_{i,j}v_{j,i}\in L_2(Q_+)$$
and thus
$$\nabla^2q\in L_2(Q_+).$$

Next, since $u$ is infinitely smooth and all its derivatives are bounded in space and time, after differentiation with respect to $x_k$, we find
$$\partial_tv_{,k}-\triangle v_{,k}=\nabla q_{,k}-{\rm div}\,F_{,k} +u_{,k}\cdot \nabla v+u\cdot \nabla v_{,k}
\in L_2(Q_+)$$
and therefore
$$\partial_t\nabla v, \,\,\nabla^3v\in L_2(Q_+).$$
Arguing in the same way, we find
$$\partial_t\nabla^k v, \,\,\nabla^{k+2}v,\,\,\nabla^{k+1}q\in L_2(Q_+)$$
for each $k=0,1,...$.

Now, we differentiate in $t$ the pressure equation
$$\triangle\partial_t q={\rm div}\,({\rm div}\,\partial_tF-\partial_tu\cdot \nabla v-u\cdot \nabla\partial_tv)$$
and establish
$$\nabla^k\partial_t q\in L_2(Q_+)$$
for any $k=1,2,...$.
Then
$$\partial^2_tv-\triangle\partial_tv=-{\rm div}\,\partial_tF+\nabla\partial_t q+
\partial_tu\cdot \nabla v+u\cdot \nabla\partial_tv$$
and thus
$$\nabla^k\partial^2_tv,\,\nabla^{k+2}\partial_tv\in L_2(Q_+)$$
for $k=0,1,...$. And so on. $\Box$

\setcounter{equation}{0}
\section{Appendix II}

We recall that 
$u(x,t)=\tilde u (x,-t)$ and $p(x,t)=-\tilde  p(x,-t)$ for $t>0$. Then
\begin{equation}\label{backNSS}
    -\partial_t u+ u\cdot\nabla  u-\triangle  u=-\nabla  p,\qquad {\rm div}\, u=0
\end{equation}
in $Q_+$ in the sense of distributions.

So, let $v$ be a solution to (\ref{moreperturbproblem}) and (\ref{moreinitialldata}).  Now, for a compactly supported smooth function $\psi$  in $Q_+$, integration by parts gives

$$\int\limits_{Q_+}u\cdot \psi{\rm div}\,{F}dx dt=$$
$$=\int\limits_{Q_+}u\cdot \psi\Big(-\partial_t{v}+u\cdot\nabla {v}+\triangle {v}+\nabla {q}\Big)dx dt=$$$$=\int\limits_{Q_+}\Big(u\cdot {v}\partial_t\psi-u\cdot {v}u\cdot\nabla\psi-u_i{v}_{i,j}\psi_{,j}+u_{i,j}
{v}_i\psi_{,j}-{q}u\cdot\nabla\psi\Big)dx dt+$$
$$+{v}\psi\cdot\Big(\partial_tu-u\cdot\nabla u+\triangle u\Big)dx dt=$$
$$=\int\limits_{Q_+}\Big(u\cdot {v}\partial_t\psi-u\cdot {v}u\cdot\nabla\psi-2u_i{v}_{i,j}\psi_{,j}+
(u_{i,j}{v}_i+u_i{v}_{i,j})\psi_{,j}-{q}u\cdot\nabla\psi\Big)dx dt+$$$$+\int\limits_{Q_-}{v}\psi\cdot \nabla pdx dt=$$
$$=\int\limits_{Q_+}\Big(u\cdot {v}\partial_t\psi-u\cdot {v}u\cdot\nabla\psi-2u_i{v}_{i,j}\psi_{,j}-u\cdot {v} \triangle\psi
-({q}u+p{v})\cdot\nabla\psi\Big)dx dt.$$


As it has been shown in \cite{SZ2006} and \cite{S2006}, one may assume that some scaled invariant energy quantities of $w$  are bounded. The same quantities remain to be bounded for $\tilde u$ and therefore for $u$. To be precise, we have 
\begin{equation}\label{scale invariant quantity}
	A+E+C+D+C_1+D_1+F+H+G=M<\infty,\end{equation}
where
$$A=\sup\limits_{R>0}\sup\limits_{R^2>t>0}\frac 1R\int\limits_{B(R)}|u(x,t)|^2dx,$$
$$E=\sup\limits_{R>0}\frac 1R\int\limits_{Q_+(R)}|\nabla u(x,t)|^2dx dt,$$
$$C=\sup\limits_{R>0}\frac 1{R^2}\int\limits_{Q_+(R)}|u|^3dx dt,\qquad D=\sup\limits_{R>0}\frac 1{R^2}\int\limits_{Q_+(R)}|p|^\frac 32dx dt,$$
$$C_1=\sup\limits_{R>0}\frac 1{R^\frac 53}\int\limits_{Q_+(R)}|u|^\frac {10}3dx dt,\qquad D_1=\sup\limits_{R>0}\frac 1{R^\frac 53}\int\limits_{Q_+(R)}|p|^\frac 53dx dt,$$
$$F=\sup\limits_{R>0}\frac 1{R^3}\int\limits_{Q_+(R)}|u|^2dx dt,\qquad H=\sup\limits_{R>0}\frac 1{R^\frac 52}\int\limits_{Q_+(R)}|u|^\frac 52dx dt,$$$$ G=\sup\limits_{R>0}\frac 1{R}\int\limits_{Q_+(R)}|u|^4dx dt$$
and $Q_+(R):=B(R)\times ]0,R^2[$.

We pick $\psi(x,t)=\chi(t)\varphi(x)$. Using simple arguments and smoothness of $u$ and $v$, we can get rid of $\chi$ and have
$$J_R(T)=\int\limits^T_0\int\limits_{\mathbb R^3}u\cdot \varphi{\rm div}\,{F}dx dt=-\int\limits_{\mathbb R^3}\varphi(x)u(x,T)\cdot {v}(x,T)dx+$$$$+\int\limits^T_0\int\limits_{\mathbb R^3}\Big(u\cdot {v}u\cdot\nabla\varphi+2u_i{v}_{i,j}\varphi_{,j}+u\cdot {v} \triangle\varphi
+({q}u+p{v})\cdot\nabla\varphi\Big)dx dt.$$

Fix a cut-off function
$\varphi(x)=\xi(x/R)$, where $\xi\in C^\infty_0(\mathbb R^3)$ with the following properties: $0\leq \xi\leq 1$, $\xi(x)=1$ if $|x|\leq 1$, and $\xi(x)=0$ if $|x|\geq 2$.
Our aim is to show that
$$J^1_R(T)=\int\limits^T_0\int\limits_{\mathbb R^3}\Big(u\cdot {v}u\cdot\nabla\varphi+2u_i{v}_{i,j}\varphi_{,j}+u\cdot {v} \triangle\varphi
+({q}u+p{v})\cdot\nabla\varphi\Big)dx dt$$
tends to zero if $R\to \infty$.

Assuming $R^2>T$, we start with
$$\Big|\int\limits^T_0\int\limits_{\mathbb R^3}2u_i{v}_{i,j}\varphi_{,j}dx dt\Big|\leq \frac cR\Big(\int\limits^T_0\int\limits_{B(2R)}|u|^2dx dt\Big)^\frac 12\Big(\int\limits^T_0\int\limits_{\mathbb R^3}|\nabla {v}|^2dx dt\Big)^\frac 12\leq$$$$\leq c\sqrt A  \sqrt {\frac TR}\|\nabla {v}\|_{2,Q_+}\to 0$$
as $R\to \infty$.

Next, we have
$$\Big|\int\limits^T_0\int\limits_{\mathbb R^3}u\cdot {v} \triangle\varphi dx dt \Big|\leq \frac c{R^2}\Big(\int\limits_T^0\int\limits_{B(2R)}|u|^2dx dt\Big)^\frac 12\Big(\int\limits^T_0\int\limits_{\mathbb R^3}| {v}|^2dx dt\Big)^\frac 12\leq $$$$\leq c\sqrt A  \sqrt \frac {T^2}{R^3}\|{v}\|_{2,\infty,Q_+}\to 0$$
as $R\to \infty$.

The third term is estimated as follows:
$$\Big|\int\limits^T_0\int\limits_{\mathbb R^3}\Big(u\cdot {v}u\cdot\nabla\varphi dx dt\Big|\leq $$$$\leq \frac c{R}\Big(\int\limits^T_0\int\limits_{B(2R)}|w|^4dx dt\Big)^\frac 12\Big(\int\limits^T_0\int\limits_{B(2R)\setminus B(R)}| {v}|^2dx dt\Big)^\frac 12\leq $$$$\leq \frac c{\sqrt R}\Big(\frac 1{2R}\int\limits_{Q_+(2R)}|u|^4dx dt\Big)^\frac 12\Big(\int\limits^T_0\int\limits_{\mathbb R^3}| {v}|^2dx dt\Big)^\frac 12\leq $$
$$\leq c{\sqrt {\frac {GT}R}}\|{v}\|_{2,\infty,Q_+}\to 0$$
as $R\to \infty$.

Now, we are going to estimate terms with pressure
$$\Big|\int\limits^T_0\int\limits_{\mathbb R^3}p{v}\cdot\nabla\varphi dx dt\Big|\leq \frac C{R}\Big(\int\limits^T_0\int\limits_{B(2R)}|p|^\frac 53dx dt\Big)^\frac 35\Big(\int\limits^T_0\int\limits_{B(2R)\setminus B(R)}| {v}|^\frac 52dx dt\Big)^\frac 25\leq$$
$$\leq cD^\frac 35\Big(\int\limits^T_0\int\limits_{B(2R)\setminus B(R)}| {v}|^\frac 52dx dt\Big)^\frac 25\to 0$$
as $R\to \infty$. The latter is true since the integral
$$\int\limits^T_0\int\limits_{\mathbb R^3}| {v}|^\frac 52dx dt$$
is finite. Indeed, this follows from the multiplicative inequality
$$\int\limits^T_0\int\limits_{\mathbb R^3}| {v}|^\frac 52dx dt \leq c T^\frac 58\|{v}\|_{2,\infty,Q_+}^\frac 74\|\nabla {v}\|_{2,Q_+}^\frac 34.$$

The most difficult term is the last one. To treat it, we split pressure ${q}$ into two parts
${q}=P_1+P_2$ so that
$$\triangle P_1=-{\rm div}\,{\rm div}\, {v}\otimes u$$
and
$$\triangle P_2={\rm div}\,{\rm div}\,F.$$
 As to the second part $P_2$, we know that it belongs to $L_\infty(0,T;L_2(\mathbb R^3))$. This  is an immediate consequence of the solution formula
$$P_2(x,t)=\frac 13 {\rm trace}\,F(x,t) -\int\limits_{\mathbb R^3}K(x-y):F(y,t) dy,$$
with the kernel $K(x)=\frac 1{4\pi}\nabla^2 \Big(\frac 1{|x|}\Big)$. Then, we have
$$\Big|\int\limits^T_0\int\limits_{\mathbb R^3}P_2u\cdot\nabla\varphi dx dt\Big|\leq \frac cR\Big(\int\limits_0^T\int\limits_{\mathbb R^3}|P_2|^2dx dt\Big)^\frac 12\Big(\int\limits_0^T\int\limits_{B(2R)}|u|^2dx dt\Big)^\frac 12\leq $$$$
\leq\sqrt A c\sqrt\frac {T^2}R\|P_2\|_{2,\infty,Q_+}
\to 0$$
as $R\to \infty$.

Regarding the second part, we are going to use the following decomposition:
$$P_1(x,t)=p_{1R}(x,t)+p_{2R}(x,t)+c_R(t),$$
where
$$p_{1R}(x,t)=-\frac 13 u(x,t)\cdot{v}(x,t)+\int\limits_{B(3R)}K(x-y):{v}(y,t)\otimes w(y,t)dy,$$
$$p_{2R}(x,t)=\int\limits_{\mathbb R^3\setminus B(3R)}(K(x-y)-K(-y)):{v}(y,t)\otimes u(y,t)dy,$$
and
$$c_R(t)=\int\limits_{\mathbb R^3\setminus B(3R)}K(-y):{v}(y,t)\otimes w(y,t)dy.$$
First of all, we observe that
$$\int\limits^T_0\int\limits_{\mathbb R^3}P_1u\cdot\nabla\varphi dx dt=
\int\limits^T_0\int\limits_{\mathbb R^3}p_{1R}u\cdot\nabla\varphi dx dt+\int\limits^T_0\int\limits_{\mathbb R^3}p_{2R}u\cdot\nabla\varphi dx dt.$$
By the theory of singular integrals,
$$\int\limits_{B(3R)}|p_{1R}|^\frac 43dx\leq c\int\limits_{B(3R)}|u|^\frac 43 |{v}|^\frac 43 dx $$
and thus
$$\int\limits_0^T\int\limits_{B(3R)}|p_{1R}|^\frac 43dx dt\leq c\Big(\int\limits_0^T\int\limits_{B(3R)}|u|^4dx dt\Big)^\frac 13 \Big(\int\limits_0^T\int\limits_{B(3R)}|{v}|^2 dx dt\Big)^\frac 23\leq $$$$
\leq cR^\frac 13G^\frac 13T^\frac 23\|{v}\|^\frac 43_{2,\infty,Q_+} .$$
So,
$$\Big|\int\limits^T_0\int\limits_{\mathbb R^3}p_{1R}u\cdot\nabla\varphi dx dt\Big|\leq
\frac cR\Big(\int\limits^0_T\int\limits_{B(2R)}|p_{1R}|^\frac 43dx dt\Big)^\frac 34
\Big(\int\limits^0_T\int\limits_{B(3R)}|u|^4dx dt\Big)^\frac 14\leq
$$
$$\leq \frac cR R^\frac 14G^\frac 14T^\frac 12 \|{v}\|_{2,\infty,Q_+}R^\frac 14G^\frac 14\to0$$
as $R\to \infty$.

Assuming that $R<|x|<2R$ and $0<t<T$, we have for the second counterpart the following estimate
$$|p_{2R}(x,t)|\leq c\int\limits_{\mathbb R^3\setminus B(3R)}\frac {|x|}{|y|^4}|u(y,t)||{v}(y,t)|dy\leq $$
$$\leq cR\sum\limits^\infty_{k=0}\frac 1{(R2^{k})^4}\int\limits_{R2^k<|y|<R2^{k+1}}|u(y,t)||{v}(y,t)|dy\leq$$
$$\leq cR\sum\limits^\infty_{k=0}\frac 1{(R2^{k})^4}\Big(\int\limits_{B(R2^{k+1})}|u(y,t)|^2dy\Big)^\frac 12\Big(\int\limits_{B(R2^{k+1})}|{v}(y,t)|^2dy\Big)^\frac 12\leq$$
$$\leq cR\Big(\int\limits_{\mathbb R^3}|v(y,t)|^2dy\Big)^\frac 12
\sum\limits^\infty_{k=0}\frac 1{(R2^{k})^4}(R2^{k+1})^\frac 12A^\frac 12\leq$$
$$\leq  \sqrt A \frac c{R^\frac 52} \|{v}\|_{2,\infty,Q_+}.$$
Then,
$$\Big|\int\limits^T_0\int\limits_{\mathbb R^3}p_{2R}w\cdot\nabla\varphi dx dt\Big|\leq
\frac cR\int\limits_0^T\sqrt A \frac 1{R^\frac 52}\|{v}\|_{2,\infty,Q_+}\int\limits_{B(2R)}|u(x,t)|dxdt\leq $$$$\leq \sqrt A\frac c{R^\frac 72}|B(2R)|^\frac 12\|{v}\|_{2,\infty,Q_+}\int\limits_0^T
\Big(\int\limits_{B(2R)}|u(y,t)|^2dy\Big)^\frac 12dt$$
$$\leq (-AT)\frac c{R^\frac 32}\|{v}\|_{2,\infty,Q_+}\to 0$$
as $R\to\infty$.
So, finally, we have
$$\int\limits_0^T\int\limits_{\mathbb R^3} u\cdot {\rm div}\,F dx dt=
-\lim\limits_{R\to\infty}\int\limits_{\mathbb R^3}\varphi(x)u(x,T)\cdot {v}(x,T)dx.$$
Taking into account $u(\cdot,T)\cdot v(\cdot,T)\in L_1(\mathbb R^3)$, see (\ref{L1}), we conclude that 
$$\int\limits_0^T\int\limits_{\mathbb R^3} u\cdot {\rm div}\,F dx dt=
-
\int\limits_{\mathbb R^3}u(x,T)\cdot {v}(x,T)dx.$$
for any $T>0$.




\end{document}